\documentclass[12pt]{article} 
\usepackage{amsmath}
\usepackage{amsthm}
\usepackage{amssymb}

\numberwithin{equation}{section}

\theoremstyle{plain}

\newtheorem{theorem}{Theorem}[section]
\newtheorem{proposition}[theorem]{Proposition}
\newtheorem{definition}[theorem]{Definition}
\newtheorem{algorithm}[theorem]{Algorithm}
\newtheorem{corollary}[theorem]{Corollary}

\newtheorem{lemma}[theorem]{Lemma}

\newcommand{\Fa}{\mathfrak{a}}

\newcommand{\Fq}{\mathfrak{q}}
\newcommand{\N}{\mathfrak{N}}
\title{N-systems, class polynomials for double $\eta$-quotients and singular values of $J$-invariant function}
\author{Shunsuke Yoshimura,  Aya Comuta  and  Noburo Ishii}
\date{ }
\begin{document}
\maketitle

\section{Introduction}\label{section1}

In elliptic curve cryptosystems, it is important to construct elliptic curves with a required number of points over finite fields. The theory of complex multiplication leads an approach to the construction of elliptic curves over finite fields. When we construct elliptic curves, we often make use of the class polynomial of $J$-invariant function. However, this method has the problem on practical use that the coefficients of the polynomial grow very large.
To overcome the problem, the polynomials with small coefficients have been invented by using Weber functions.
Recently, Enge and Schertz \cite{EN1} gave the method of using the double $\eta$-quotient.
The double $\eta$-quotient is not $\rm{SL}_{2}(\mathbb{Z})$-invariant but $\Gamma^{0}(N)$-invariant. For this reason, they considered $N$-systems and calculated  class polynomials with respect to double $\eta$-quotients.
Furthermore, they can calculate the modular invariants  $J(\Fa)$, for an ideal $\Fa$ of an imaginary quadratic field, by using the modular polynomial from the class polynomials with respect to double $\eta$-quotients. However in their method, it is necessary to count the number of rational points of elliptic curves corresponding to solutions of the modular equation over a finite field, because in advance we can not know which solution of the modular equation is that corresponding to the modular invariant.
In this article, we shall give a method to reduce the amount of computation in the process of counting the number of rational points. Thus, if we are in the situation that the modular invariant is a multiple root of the modular polynomial, we may expect that the amount of computation reducers to at most a half of the original one. 
In Section 2, we give some basic results and definitions. In Section 3, we determine the relations between class polynomials with respect to double $\eta$-quotients and $N$-systems.
In Section 4 , we give a condition that the modular invariant is a multiple root of the modular polynomial. In section 5, we give an example.

\section{Some basic results and definitions}\label{section2}

For two prime numbers $p_1$ and $p_2$, the double $\eta$-quotient $\mathfrak{w}_{p_{1},p_{2}}(z)$ of level $N=p_1p_2$
is defined by 
$$
 \mathfrak{w}_{p_{1},p_{2}}(z)  =  \frac{\eta(z/p_1) \eta(z/p_2)}{\eta(z) \eta(z/p_{1}p_{2})},
$$
where $\eta(z)$ denotes the Dedekind $\eta$-function defined by
$$
 \eta(z) = q^{1/24} \prod_{n=1}^{\infty}(1-q^{n}), \qquad 
 q = q(z) = e^{2\pi iz}.
$$
The function $\mathfrak{w}_{p_{1},p_{2}}(z)$ is invariant under the modular group 
$$\Gamma^{0}(N) = \left\{ 
  \left(\left.
   \begin{array}{cc}
    a & b \\
    c & d \\
   \end{array}
  \right)
  \in \rm{SL}_{2}(\mathbb{Z})\right | b\equiv 0\pmod N\right\}$$
(see \cite{MN}).
For a divisor $Q$ of $N$ such that $(Q,\frac{N}{Q})=1$,
set
$$
W_Q = 
\left(
 \begin{smallmatrix}
    -Q & N \\
    -y & -Qx \\
 \end{smallmatrix}
\right),
$$
where $x,y \in \mathbb{Z}$ and $\det (W_Q)=Q$. Then we know $W_Q$ is a normalizer of $\Gamma^{0}(N)$. Especially the normalizer $\displaystyle W_N=\begin{pmatrix}0&N\\-1&0\end{pmatrix}$ is called the Atkin-Lehner involution of $\Gamma^{0}(N)$. By Theorem 2 of \cite{EN2}, we know 
\begin{equation}\label{eq1}
 \mathfrak{w}_{p_{1},p_{2}}(W_N(z))  =  \mathfrak{w}_{p_{1},p_{2}}(z).
\end{equation}
For $\mathfrak{w}_{p_{1},p_{2}}^{s}(z)$ with $s = 24/\gcd(24,(p_1-1)(p_2-1))$, 
Enge and Schertz \cite{EN1} defined the class polynomials $H_\N(X)$ associated with $N$-systems $\N$. In the following, we shall recall their results. Let $\mathcal{O}_f$ be the order of conductor $f$ in a imaginary quadratic field $K=\mathbb{Q}(\sqrt{m})$. Let $d_K$ be a discriminant of $K$. Then the discriminant $D$ of $\mathcal{O}_f$ is given by $D=f^2d_k$. 
Let $\mathfrak{H}_f$ be the (proper) ideal class group of $\mathcal{O}_f$.
To a proper ideal $\mathfrak{a}=[\beta_1,\beta_2] = \mathbb{Z}\beta_1 + \mathbb{Z}\beta_2$ of $\mathcal{O}_f$, we associate its basis quotient
$$\alpha_{\mathfrak{a}}=\displaystyle{\frac{\beta_1}{\beta_2}} \qquad \left(\alpha_{\mathfrak{a}} \in \mathbb{H} = \{ \tau \in \mathbb{C}~|~ Im(\tau) > 0 \}\right),$$
and a quadratic form $\mathfrak q_{\Fa}(X,Y)=N_{K/\mathbb Q}(\beta_1X+\beta_2Y)/N_{K/\mathbb Q}(\mathfrak{a})$. It is noted that the basis quotient $\alpha_\Fa$ is determined up to $\rm{SL}_{2}(\mathbb{Z})$-equivalence and the form $\Fq_{\Fa}(X,Y)$ is a primitive quadratic form with integral coefficients of discriminant $D$. The following results are well known (see \cite{AC}).

Let $\Fa_1$ and $\Fa_2$ are proper ideals of $O_f$. We write $ \Fa_1 \sim \Fa_2$ if  $\Fa_1$ and $\Fa_2$ are in the same ideal class of $\mathcal{O}_f$.
Then
\begin{eqnarray}
  \Fa_1 \sim \Fa_2 
 & \Longleftrightarrow & \alpha_{\Fa_1} \; \text{is} \;\rm{SL}_{2}(\mathbb{Z})\text{-equivalent} \; \text{to} \;  \alpha_{\Fa_2} \nonumber\\
 & \Longleftrightarrow & \mathfrak q_{\Fa_1}(X,Y) \; \text{is proper equivalent to} \; \mathfrak q_{\Fa_2}(X,Y).\label{equiv1}
\end{eqnarray}
Further the map $\Fa\mapsto \mathfrak q_{\Fa}(X,Y)$ gives rise to a bijection between  $\mathfrak{H}_f$ and the proper equivalent classes of quadratic forms of discriminant $D$. Since for an ideal $\Fa$  the value $J(\alpha_\Fa)$ is independent on the choice of basis quotients $\alpha_\Fa$, we shall denote by $J(\Fa)$ the value $J(\alpha_\Fa)$. Hereafter for a triple $(A,B,C)$ of integers such that $A\ne 0$, we shall denote by $[A,B,C]$ a quadratic form $AX^2+BXY+CY^2$. Furthermore, for a quadratic form $\mathfrak q=[A,B,C]$, we put $\Fa_\Fq=[A,\frac{-B+\sqrt D}2]$ and $\alpha_\Fq=\frac{-B+\sqrt D}{2A}$. Let $h(D)$ be the class number of the order $\mathcal{O}_f$. Then we know  there exist $h(D)$ isomorphic classes of elliptic curves with complex multiplication by $\mathcal{O}_f$. They are represented  by the elliptic curves with $j$-invariants $j(\alpha_{\mathfrak{a_i}})$, where $\Fa_i~(i=1,\dots ,h(D))$ are ideals representing all classes of $\mathfrak{H}_f$. Let $K_f$ be the ring class field of $K$ of conductor $f$. Then the theory of complex multiplication shows  $J(\Fa_i)$ generates $K_f$ over $K$ for each $i$ and $J(\Fa_i),~(i=1,\dots ,h(D))$ are conjugate  to each other over $\mathbb Q$. To calculate the values $J(\Fa)$, we use the classical class equation for the order $\mathcal{O}_f$ defined by
\[H_D[J](X)=\prod_{i=1}^{h(D)} (X-J(\Fa_i)).\]
However, since this polynomial has large integral coefficient, it is hard to compute the polynomial for large $D$'s.  Enge and Schertz  \cite{EN1} devised the method of using the double $\eta$-quotient to obtain a class polynomial with small integral coefficients. Since double $\eta$-quotient is $\Gamma^{0}(N)$-invariant but not $\rm{SL}_{2}(\mathbb{Z})$-invariant, the values $\mathfrak{w}_{p_{1},p_{2}}^s(\alpha_\Fa)$ depend on the choice of ideals $\Fa$ in an ideal class. Therefore to define a class polynomial for $\mathfrak{w}_{p_{1},p_{2}}^s$, we must use the $N$-systems introduced by Schertz \cite{RS}.   
 
\begin{definition}
Let $\N$ be a set of $h(D)$ primitive quadratic forms $[A_i,B_i,C_i]$ of discriminant $D$. we say 
$\N$ an $N$-system for $D$ if quadratic forms $[A_i,B_i,C_i]$ satisfy the following conditions:
\begin{enumerate}
  \item $\gcd(A_i,N)=1, B_i\equiv B_j\pmod{2N},~N|C_i$ for every $i,j$,
 \item the ideals $[A_i, \frac{-B_i+\sqrt{D}}{2}]$ form a representative system of $\mathfrak{H}_f$.
\end{enumerate}
\end{definition}
The following result is basic for $N$-systems. See Theorems 3.1, 3.2 of \cite{EN1}.
\begin{theorem}\label{theorem2.3} 
Assume prime numbers $p_1$ and $p_2$ satisfy the following conditions:
\begin{enumerate}
 \item If $p_1 \neq p_2$, then $\bigl( \frac{D}{p_1} \bigr),\bigl( \frac{D}{p_2} \bigr) \neq -1$,
 \item if $p_1=p_2=p$, then either $\bigl( \frac{D}{p} \bigr)=1$ or $p|f$.
\end{enumerate}
Then there exists an $N$-system $\N=\{\Fq_i\}$ for $D$ . Further $\mathfrak{w}_{p_{1},p_{2}}^{s}(\alpha_{\Fq_i}) \in K_f$ for every $i$ and  $\mathfrak{w}_{p_{1},p_{2}}^{s}(\alpha_{\Fq_i})$ are conjugate under $G(K_f/K)$ to each other. 
\end{theorem}
Now the class polynomial of $\mathfrak{w}_{p_{1},p_{2}}^{s}$ associated with an $N$-system $\N$ is defined by 
\[
H_\N(X) = \prod_{i=1}^{h(D)} \left (X-\mathfrak{w}_{p_{1},p_{2}}^{s}(\alpha_{\Fq_i})\right).
\]
By Theorem~\ref{theorem2.3}, we have $H_\N(X)\in K[X]$. Further Corollary 3.1 of \cite{EN1} gives
\begin{corollary}\label{Cor2.4} 
 Suppose that following conditions hold:
\begin{enumerate}
 \item If $p_1 \neq p_2$, then $\bigl( \frac{D}{p_1} \bigr),\bigl( \frac{D}{p_2} \bigr) \neq -1$, and $p_1,p_2 \nmid f$;
 \item if $p_1=p_2=p \neq 2$, then $\bigl( \frac{D}{p} \bigr)=1$ or $p|f$;
 \item if $p_1=p_2=2$, then $\bigl( \frac{D}{2} \bigr)=1$, or $2|f$, but $D\not\equiv 4\pmod{32}$.
\end{enumerate}
Then $H_\N(X)\in\mathbb{Z}[X]$.
\end{corollary}
To obtain modular invariants $J(\Fa_{\Fq})$ from the singular values $\mathfrak{w}_{p_{1},p_{2}}^{s}(\alpha_{\Fq})$, we use a modular polynomial $\Phi_{p_1,p_2}$ that associates $\mathfrak{w}_{p_{1},p_{2}}$ with $J$, which is defined as follows.

\begin{definition}
$$
\Phi_{p_1,p_2}(X,J) = \prod_{\sigma \in Iso(\mathbb{C}_{\Gamma^{0}(N)}/\mathbb{C}_{\Gamma})}
(X-\sigma(\mathfrak{w}_{p_{1},p_{2}}^{s})),
$$
where $\mathbb{C}_{\Gamma^{0}(N)}$ and $\mathbb{C}_{\Gamma}$ denote the modular function field 
of $\Gamma^{0}(N)$ and $\Gamma=\rm{SL}_{2}(\mathbb{Z})$ respectively.
\end{definition}

We know that $\Phi_{p_1,p_2}(X,J) \in \mathbb{Z}[X,J]$ and 
$\Phi_{p_1,p_2}(X,J)$ is minimal polynomial of $\mathfrak{w}_{p_{1},p_{2}}^{s}$ over $\mathbb{C}(J)$
by Theorems 7 and 8 by \cite{EN2}. To obtain elliptic curves with complex multiplication by $\mathcal{O}_f$ over finite field $\mathbb{F}_q$ of $q$-elements, the polynomials $H_\N(X)$ and $\Phi_{p_1,p_2}(X,J)$ are used in the following algorithm.

Assume that $q$ is a prime number which splits completely in $K_f$.

\begin{algorithm}

\begin{enumerate}
 \item Construct an $N$-system $\{\Fq_i\}$.
 \item Compute $\mathfrak{w}_{p_{1},p_{2}}^{s}(\alpha_{\Fq_i})$ and $H_\N(X)
$.
 \item Compute a roots $\overline{\mathfrak{w}}$ of $H_\N(X) \bmod q $.
 \item Compute  $\mathbb{F}_q$-rational roots $\overline{J_{k}}$ of $\Phi_{p_1,p_2}(\overline{\mathfrak{w}},J) \bmod q $.
 \item Output the desired $J$-invariant among $\overline{J_{k}}$.
\end{enumerate}

\end{algorithm}

In step 5, it is necessary to count the number of $\mathbb{F}_q$-rational points of each elliptic curve $E_k$ with the j-invariant $\overline{J_{k}}$ to determine which elliptic curve $E_k$ has complex multiplication by $\mathcal{O}_f$. In step 4, if $\Phi_{p_1,p_2}(\overline{\mathfrak{w}},J) \bmod q $ has the degree 2 in $J$ and the multiple root $\overline{J}$ in $J$, then it is not necessary to count the number of rational points. Accordingly, we will consider the condition that the polynomial has the multiple root in Section \ref{section4}.

\section{$N$-systems and class polynomials}\label{section3}

In this section, we study the relation between the class polynomials and $N$-systems.

\
\begin{lemma}\label{lemma1}
Let $\{[A_i,B_i,C_i]\}$ and $\{[A_i',B_i',C_i']\}$ be $N$-systems for $D$.
Suppose that $[A_i,B_i,C_i]$ and $[A_i',B_i',C_i']$ are proper equivalent.
Then $\frac{-B_i+\sqrt{D}}{2A_i}$ is $\Gamma^0(N)$-equivalent to $\frac{-B_i'+\sqrt{D}}{2A_i'}$ if and only if $B_i\equiv B_i'\pmod {2N}$. In particular, if $B_i\equiv B_i'\pmod {2N}$, then 
$$
\mathfrak{w}_{p_{1},p_{2}}^{s}(\frac{-B_i+\sqrt{D}}{2A_i})=\mathfrak{w}_{p_{1},p_{2}}^{s}(\frac{-B_i'+\sqrt{D}}{2A_i'}).
$$
\end{lemma}

\begin{proof}
Since $[A_i,B_i,C_i]$ and $[A_i',B_i',C_i']$ are proper equivalent, there exists a matrix 
$
 M= \left(
   \begin{array}{cc}
    a & b \\
    c & d \\
   \end{array}
  \right)
  \in \rm{SL}_{2}(\mathbb{Z})
$
such that
\begin{eqnarray}
\frac{-B_i+\sqrt{D}}{2A_i} = \frac{a(\frac{-B_i'+\sqrt{D}}{2A_i'})+b}{c(\frac{-B_i'+\sqrt{D}}{2A_i'})+d}. \label{5}
\end{eqnarray}
We have only to prove
$M \in \Gamma^0(N)$ if and only if $B_i\equiv B_i'\pmod {2N}$.
By (\ref{5}), we have
\begin{eqnarray}
B_iB_i'c-2A_i'B_id+Dc & = &- 2A_iB_i'a+4A_iA_i'b,  \label{6}\\
-B_ic-B_i'c+2A_i'd & = & 2A_ia. \label{7}
\end{eqnarray}

By substituting (\ref{7}) into (\ref{6}), we obtain
\begin{eqnarray}
A_ia(B_i'-B_i)-2A_iC_ic=2A_iA_i'b. \label{7.1}
\end{eqnarray}

Since $\gcd(A_iA_i',N)=1$ and $N|C_i$, we have $a(B_i'-B_i)\equiv 2A_i'b\pmod {2N}$. Therefore $M \in \Gamma^0(N)$ if and only if $B_i\equiv B_i'\pmod {2N}$.
\end{proof}

The following result is deduced from Proposition 3 of \cite{RS}.
\begin{proposition}\label{prop_rs}
Let $[A,B,C]$ be a primitive quadratic form of discriminant $D$ such that $A>0, \gcd(A_,N)=1$ and $N|C$. Then there exists an $N$-system $\N$ for $D$ containing $[A,B,C]$. In particular, for an integer $B$ such that $B^2\equiv D\pmod {4N}$, there exists an $N$-system $\N$ for $D$ containing $[1,B,(B^2-D)/4]$.
\end{proposition}

By Proposition~\ref{prop_rs} and Lemma~\ref{lemma1}, we know the class polynomials of double $\eta$-quotient associated with $N$-systems depend only on integers $B$, considered mod $2N$, such that $B^2\equiv D\pmod {4N}$. Thus, hereafter, we shall denote by $H_{B,N}(X)$ the class polynomial $H_\N(X)$ associated with  an $N$-system $\N$ containing a form $[1,B,(B^2-D)/4]$. In the following, we shall fix an $N$-system containing $[1,B,(B^2-D)/4]$ and shall denote it by $\N_B$.
\begin{lemma}\label{lemma2}
Assume that $p_1$ and $p_2$ are odd primes. Let $N(D)$ be the number of integers $B \bmod 2N$ such that $B^2\equiv D \pmod {4N}$.
Then
\begin{eqnarray}
N(D) = 
\begin{cases}
4 \quad if \quad {\bigl( \frac{D}{p_1} \bigr)=1 \quad and \quad \bigl( \frac{D}{p_2} \bigr)=1}, \\
2 \quad if \quad {\bigl( \frac{D}{p_1} \bigr)=1 \quad and \quad \bigl( \frac{D}{p_2} \bigr)=0}.
\end{cases}
\end{eqnarray}
\end{lemma}

\begin{proof}
Let us consider the case  $\bigl( \frac{D}{p_1} \bigr)=1$  and  $\bigl( \frac{D}{p_2} \bigr)=1$. Then there exists an integer $\alpha_i$ such that $\alpha_i^2\equiv D \pmod{p_i}$ for $i=1,2$. By Chinese reminder theorem, we see $B^2\equiv D \pmod{4N}$ if and only if
\begin{eqnarray}
 B & \equiv & D \pmod{2}, \\
 B & \equiv & \pm \alpha_i \pmod{p_i}\quad (i=1,2). \label{1} 
\end{eqnarray}
This shows $N(D) = 4$. The remaining case can be treated similarly. Thus we omit the details.
\end{proof}
If $N$ is odd, we obtain at most $N(D)$ distinct class polynomials of the double $\eta$-quotient associated with $N$-systems for $D$.
\begin{lemma}\label{lemma3}
Let $B$ be an integer such that $B^2\equiv D \pmod {4N}$. Then $H_{-B,N}(X)=H_{B,N}(X)$.
\end{lemma}

\begin{proof}
We know the $q$-expansion of $\mathfrak{w}_{p_{1},p_{2}}^{s}(z)$ is rational, thus
$$
\mathfrak{w}_{p_{1},p_{2}}^{s}(z) = \sum a_n q^n \qquad (a_n \in \mathbb{Q}, \quad q=e^{2\pi iz})
$$
(see section 3 of \cite{EN2}). Therefore we have  
$$ 
\overline{\mathfrak{w}_{p_{1},p_{2}}^{s}(z)} = \sum a_n \overline{q}^n = \mathfrak{w}_{p_{1},p_{2}}^{s}({\overline{q}})=\mathfrak{w}_{p_{1},p_{2}}^{s}(-\overline z)
$$
and
\[
\overline{\mathfrak{w}_{p_{1},p_{2}}^{s}\Bigl(\frac{-B+\sqrt{D}}{2}\Bigr)} = \mathfrak{w}_{p_{1},p_{2}}^{s}\Bigl(\frac{B+\sqrt{D}}{2}\Bigr).
\]
Since  Corollary~\ref{Cor2.4} shows that $H_{B,D}(X)\in\mathbb Z[X]$, we have $H_{B,D}(X)=\overline{H_{B,N}(X)}=H_{-B,N}(X)$.

\end{proof}

\begin{theorem}\label{theorem2}
Let $N(H_D)$ be the number of distinct class polynomials $H_{B,N}(X)$ associated with $N$-systems.
Then
\begin{eqnarray}
N(H_D) = 
\begin{cases}
1,2 \quad &if \quad {\bigl( \frac{D}{p_1} \bigr)=1 \quad and \quad \bigl( \frac{D}{p_2} \bigr)=1}, \\
1 \quad &if \quad {\bigl( \frac{D}{p_1} \bigr)=1 \quad and \quad \bigl( \frac{D}{p_2} \bigr)=0}.
\end{cases}
\end{eqnarray}
\end{theorem}

\begin{proof}
Lemmas \ref{lemma2} and \ref{lemma3} imply that $N(H_D) \le N(D)/2$. Thus we have the assertion.
\end{proof}

In the case $N(H_D)=2$, we have two class polynomials $H_{B,N}(X)$ and $H_{B',N}(X)$,
where $B,B'$ are integers such that $B^2\equiv D\pmod{4N},~B'\equiv B\pmod{p_1},~B'\equiv -B\pmod{p_2}$. We shall show $H_{B',N}(X)$ is obtainable from $H_{B,N}(X)$ by a simple transformation. 
We shall use the following transformation formula of the Dedekind $\eta$-function (see Theorem 1 of \cite{EN2}). 
\begin{theorem}\label{EN2_Th1}  
Let 
$
M = 
  \left(
   \begin{array}{cc}
    a & b \\
    c & d \\
   \end{array}
  \right)
  \in \rm{SL}_{2}(\mathbb{Z})
$
be normalized such that $c\geq 0$, and $d>0$ if $c=0$.
Write $c=\gamma 2^{\lambda}$ with $\gamma$ odd; by convention, $\gamma=\lambda=1$ if $c=0$. Then
$$
\eta(Mz)  =  \epsilon(M)\sqrt{cz+d} \eta(z) 
$$
with 
$$
\mathfrak{R}(\sqrt{cz+d})>0 , \quad
\epsilon(M)  = \Bigl( \frac{a}{\gamma} \Bigr) \zeta^{ab+c(d(1-a^2)-a)+3\gamma(a-1)+\frac{3}{2}\lambda(a^2-1)}_{24}.
$$
\end{theorem}

\begin{lemma}\label{lemma_involution}
Let $
W_{p_1} = 
  \left(
   \begin{array}{cc}
    -p_1 & N \\
    -y & -p_1x \\
   \end{array}
  \right)
$
with $y < 0,~y\equiv 1\pmod 2$ and $p_1x+p_2y=1$. Then
$$
 \mathfrak{w}_{p_1,p_2}^s(W_{p_1}(z)) = \frac{\Bigl( \frac{p_1}{p_2} \Bigr)^s}{\mathfrak{w}_{p_1,p_2}^{s}(z)}.
$$
\end{lemma}

\begin{proof}

 By Theorem \ref{EN2_Th1},
\begin{eqnarray*}
 \mathfrak{w}_{p_1,p_2}(W_{p_1}(z)) = \mathfrak{w}_{p_1,p_2}(\frac{-p_1z+N}{-yz-p_1x})
 = \frac{\eta(\frac{-p_1z+N}{p_1(-yz-p_1x)})\eta(\frac{-p_1z+N}{p_2(-yz-p_1x)})}{\eta(\frac{-p_1z+N}{-yz-p_1x})\eta(\frac{-p_1z+N}{p_1p_2(-yz-p_1x)})}
 = \frac{\epsilon^{*}}{\mathfrak{w}_{p_1,p_2}(z)}.
\end{eqnarray*}
Here
$$
\epsilon^{*} ={\footnotesize \frac{
\epsilon
  \left(\Bigl(
   \begin{array}{cc}
    -1 & p_2 \\
    -y & -p_1x \\
   \end{array}\Bigr)
  \right)
\epsilon
  \left(\Bigl(
   \begin{array}{cc}
    -p_1 & 1 \\
    -p_2y & -x \\
   \end{array}\Bigr)
  \right)
}{
\epsilon
  \left(\Bigl(
   \begin{array}{cc}
    -p_1 & p_2 \\
    -y & -x \\
   \end{array}\Bigr)
  \right)
\epsilon
  \left(\Bigl(
   \begin{array}{cc}
    -1 & 1 \\
    -p_2y & -p_1x \\
   \end{array}\Bigr)
  \right)
}}
= \frac{\Bigl(\frac{-1}{-y} \Bigr)\zeta^{a}_{24} \Bigl( \frac{-p_1}{-p_2y} \Bigr)\zeta^{b}_{24}}
{\Bigl(\frac{-p_1}{-y} \Bigr)\zeta^{c}_{24} \Bigl( \frac{-1}{-p_2y} \Bigr)\zeta^{d}_{24}}
= \Bigl(\frac{p_1}{p_2} \Bigr)\zeta^{a+b-c-d}_{24}
$$
with
\begin{eqnarray*}
a & = & -p_2-y+3\gamma(-2), \\
b & = & -p_1 -p_2y(-x(1-p_1^2)+p_1)+3p_2\gamma(-p_1-1)+\frac{3}{2}\lambda(p_1^2-1), \\
c & = & -p_1p_2 -y(-x(1-p_1^2)+p_1)+3\gamma(-p_1-1)+\frac{3}{2}\lambda(p_1^2-1), \\
d & = & -1-p_2y+3p_2\gamma(-2).
\end{eqnarray*}
Therefore we see
$$
a+b-c-d = (1-p_2)(1-p_1)(1-y-xy(p_1+1)-3\gamma).
$$
Since $s(p_1-1)(p_2-1) \equiv 0 \pmod{24}$, we have  $(\epsilon^{*})^s= \Bigl(\frac{p_1}{p_2} \Bigr)^s$.

\end{proof}

\begin{proposition}\label{prop5}
Suppose that $\bigl( \frac{D}{p_1} \bigr)=\bigl( \frac{D}{p_2} \bigr)=1$. 
Then
$$
H_{B',N}(X)= \frac{X^{h(D)}}{H_{B,D}(0)} H_{B,N}\Bigl( \frac{\Bigl( \frac{p_1}{p_2} \Bigr)^s}{X} \Bigr).
$$
\end{proposition}

\begin{proof}
 Let $\Fq=[A,B,C]$ be a form  of the $N$-system $\N_B$. Put $\alpha_\Fq=\frac{-B+\sqrt D}{2A}$. We can write $W_{p_1}(\alpha_\Fq)=\frac{-B'+\sqrt d}{2A'}$. We see easily that 
$B'\equiv B\pmod{p_1}$, $B'\equiv -B\pmod{p_2}$ and $C'=(B'^2-D)/4A'\equiv 0\pmod N$. If $\gcd(A',N)>1$, then we shall show there exists an element $\gamma \in \Gamma^{0}(N)$ such that the first coefficient $A''$ of the quadratic form $[A'',B'',C'']$
corresponding to $\gamma W_{p_1}(\alpha)$ is prime to $N$. It is noted by the proof of Lemma~\ref{lemma1} that $B'\equiv B''\pmod{2N}$ and $N|C''$. Since $\mathfrak{w}_{p_{1},p_{2}}(\gamma W_{p_1}(\alpha_\Fq)) = \mathfrak{w}_{p_{1},p_{2}}(W_{p_1}(\alpha_\Fq))$, we have 
\[
H_{B',N}(X)=\prod_{\Fq\in\N_B}(X- \mathfrak{w}_{p_{1},p_{2}}(W_{p_1}(\alpha_\Fq))).
\]
By Lemma~\ref{lemma_involution}, we have our result. Thus we have only to prove the existence of the above $\gamma$ in the case $A'$ is not prime to $N$. Set 
$\gamma = 
  \left(
   \begin{array}{cc}
    r & N \\
    t & 1 \\
   \end{array}
  \right)
$, for integers $r,t$. Then 
\begin{eqnarray*}
\gamma W_{p_1}(\alpha) & = & \frac{r \Bigl(\frac{-B'+\sqrt{D}}{2A'}\Bigr)+N}{t\Bigl( \frac{-B'+\sqrt{D}}{2A'}\Bigr) + 1}
= \frac{(-rB'-tB'N+2A'N+2rtC')+\sqrt{D}}{2(-tB'+A'+t^2C')}.
\end{eqnarray*}
Therefore we have $A''= -tB'+A'+t^2C'$. Since $N|C'$, we know $\gcd(A'',N)=1$ if and only if $\gcd(-tB'+A',N)=1$. Assume $N|A'$. Then $\gcd(D,N)=1$ implies $(B',N)=1$. Therefore we can take $r=N+1$, $t=1$. Next assume $p_i|A'$ and $p_j\nmid A'$. Then we have $p_i\nmid B'$. Therefore we can take $r=p_jN+1$, $t=p_j$. This completes our proof.
\end{proof}

\section{Multiple roots of modular equation }\label{section4}
In this section, we assume the conditions in Corollary~\ref{Cor2.4}.
We shall study a condition that for singular values $\alpha$ of double $\eta$-quotients the polynomial $\Phi_{p_1,p_2}(\alpha,J)$ of $J$ has a multiple root.
\begin{proposition}\label{prop1}
The polynomial $\Phi_{p_1,p_2}(X,J)$ has degree $\displaystyle{N \prod_{p|N} \bigl(1+\frac{1}{p}\bigr)}$ as a polynomial of $J$ and has degree $\displaystyle{\frac{s(p_1-1)(p_2-1)}{12}}$ as a polynomial of $X$. For $\tau\in\mathbb H$, the equation $\Phi_{p_1,p_2}(\mathfrak{w}_{p_1,p_2}^s(\tau),J)=0$ has two roots $J(\tau)$ and $J(W_N (\tau))$.
In particular, if $J(\tau)=J(W_N (\tau))$, then the equation has a multiple root.
\end{proposition}

\begin{proof}
The assertion concerning the degree follows from Theorem 9 of \cite{EN2}. The equation $\Phi_{p_1,p_2}(\mathfrak{w}_{p_1,p_2}^s(\tau),J)=0$ obviously has the root $J(\tau)$.
Similarly, $J(W_N(\tau))$ is a root of $\Phi_{p_1,p_2}(\mathfrak{w}_{p_1,p_2}^s(W_N(\tau)),J)=0$. By \eqref{eq1}, we have $\mathfrak{w}_{p_{1},p_{2}}^{s}(W_N(\tau))=\mathfrak{w}_{p_{1},p_{2}}^{s}(\tau)$. Therefore $J(W_N (\tau))$ is a root of $\Phi_{p_1,p_2}(\mathfrak{w}_{p_1,p_2}^s(\tau),J)=0$.
\end{proof}

Let $\Fq=[A,B,C]$ be a form of an $N$-system $\N_B$. Since
\begin{eqnarray*}
W_N(\alpha_\Fq) & = & \frac{N}{-\frac{-B+\sqrt{D}}{2A}}
               = \frac{2AN(B+\sqrt{D})}{B^2-D}
               = \frac{B+\sqrt{D}}{2(\frac{C}{N})},
\end{eqnarray*}
the action of $W_N$ on the ideal $\mathfrak{a_\Fq}$ is given by $W_N(\mathfrak{a_\Fq}) = [\frac{C}{N}, \frac{B+\sqrt{D}}{2}]$.

\begin{lemma}\label{lemma_a_b}
If we set $\mathfrak{a}_B=[N,\frac{-B+\sqrt{D}}{2}]$, then $W_N(\mathfrak{a_\Fq}) \sim \mathfrak{a_\Fq}\mathfrak{a}_B$.
\end{lemma}

\begin{proof}Since $(A,B,C)=1$,
\begin{eqnarray*}
\overline{\mathfrak{a_\Fq}}W_N(\mathfrak{a_\Fq}) 
& = & [A,\frac{B+\sqrt{D}}{2}][\frac{C}{N},\frac{B+\sqrt{D}}{2}] \\
& = & [\frac{AC}{N},\frac{C}{N}(\frac{B+\sqrt{D}}{2}),A(\frac{B+\sqrt{D}}{2}),B(\frac{B+\sqrt{D}}{2})] \\
& = & [\frac{AC}{N},\frac{B+\sqrt{D}}{2}] \\
& = & \frac{1}{N}(\frac{B+\sqrt{D}}{2})[N,\frac{-B+\sqrt{D}}{2}].
\end{eqnarray*}
Thus we have $\overline{\mathfrak{a_\Fq}}W_N(\mathfrak{a_\Fq})\sim\mathfrak{a_\Fq}\mathfrak{a}_B$. Since $\overline{\mathfrak{a_\Fq}}\mathfrak{a_\Fq}\sim 1$, this proves the assertion.
\end{proof}

\begin{proposition}\label{prop2}
 Let $\Fq=[A,B,C]$ be a form of an $N$-system $\N_B$.  Then
$J(W_N(\alpha_\Fq))$ $=J(\alpha_\Fq)$ if and only if there exist $u,v \in \mathbb{Z}$ such that 
\begin{eqnarray}\label{con1}
\begin{cases}
u^2 - Dv^2 = 4N \\
u-Bv \equiv 0 \pmod{2N}. 
\end{cases}
\end{eqnarray}
\end{proposition}

\begin{proof}
By \eqref{equiv1} and Lemma~\ref{lemma_a_b}, we have
\begin{eqnarray*}
J(W_N(\Fa_\Fq))=J(\Fa_\Fq)
\quad & \Leftrightarrow & \quad W_N(\mathfrak{a_\Fq}) \sim \mathfrak{a_\Fq} \\
\quad & \Leftrightarrow & \quad \mathfrak{a}_B \sim 1.
\end{eqnarray*}
Further we know
the condition $\mathfrak{a}_B \sim 1$ is equivalent to the existence of an element
$
  \left(
   \begin{array}{cc}
    x & y \\
    z & w \\
   \end{array}
  \right)
  \in \rm{SL}_{2}(\mathbb{Z})
$
such that
\begin{eqnarray}\label{con2}
\frac{-B+\sqrt{D}}{2N} = \frac{x(\frac{-B+\sqrt{D}}{2})+y}{z(\frac{-B+\sqrt{D}}{2})+w}.
\end{eqnarray}
Let us assume  \eqref{con2}. Then  we have
\begin{eqnarray}
zB^2+zD-2wB & = & -2xNB+4Ny, \label{8} \\
w & = & xN+zB. \label{9}
\end{eqnarray}
By substituting (\ref{9}) into (\ref{8}), we have $y = -A\frac{C}{N}z$.
Therefore, from $xw - yz = 1$, we obtain $(Bz+2xN)^2-Dz^2=4N$.
Now, we put $u=Bz+2Nx$, $v=z$. Then we have
$u^2-Dv^2=4N$ and $x=\frac{u-Bv}{2N}$. Further, since $x \in \mathbb{Z}$, we have $u-Bv \equiv 0 \pmod{2N}$. Conversely, let $u,v$ be integers satisfying \eqref{con1}. Put $x=\frac{u-Bv}{2N}$, $y= -A\frac{C}{N}v$, $z=v$ and $w=xN+zB$. Then we have $xw-zy=1$ and \eqref{con2}. 

\end{proof}
Immediately from  Proposition~\ref{prop2}, we obtain 
\begin{corollary}\label{coro4}
If $J(W_N(\Fa_\Fq))=J(\Fa_\Fq)$, then $D > -4N$.
\end{corollary}

\begin{proof}
The condition \eqref{con1} shows
$
D = \frac{u^2-4N}{v^2} > -4N.
$
\end{proof}
\begin{proposition}\label{prop3}
Assume that there exist integers $u$ and $v$ satisfying \eqref{con1}. 
Then the equation $\Phi_{p_1,p_2}(\mathfrak{w}_{p_{1},p_{2}}^{s}(\alpha_\Fq),J)=0$ has
a multiple root $J(\Fa_\Fq)$.
\end{proposition}

\begin{proof}
The assertions is obvious.
\end{proof}

By a similar argument in Proposition~\ref{prop2}, we have:

\begin{proposition}\label{prop5}
Let the notation be as in Proposition~\ref{prop2}. Then
 ${W_N}^2(\mathfrak{a}) \sim \mathfrak{a} $ if and only if there exist integers $X,Y$ such that $Y\neq 0$ and 
\begin{eqnarray*}
\begin{cases}
X^2-DY^2=4N^2 \\
X-BY \equiv 0 \pmod{2N} \\
 (\frac{X-BY}{2N})^2 \equiv 1 \pmod{Y}.
\end{cases}
\end{eqnarray*}
\end{proposition}

\begin{theorem}\label{theorem4.7}
Let $\mathfrak{a_{\Fq_i}}~(i=1,\dots,h(D))$ be the ideals associated with the quadratic forms $\Fq_i=[A_i,B_i,C_i]$ of $\N_B$.
Then 
\begin{eqnarray*}
J(\mathfrak{a_{\Fq_1}})=J(W_N(\mathfrak{a_{\Fq_1}})) \quad \Leftrightarrow \quad J(\mathfrak{a_{\Fq_i}})=J(W_N(\mathfrak{a_{\Fq_i}}))
\quad (i=1,\dots,h(D)).
\end{eqnarray*}
\end{theorem}

\begin{proof}
We set $\mathfrak{a}_{B_i}=[N,\frac{-B_i+\sqrt{D}}{2}]$.
Since $B_1  \equiv B_i \pmod{2N}$, we know
$$
\mathfrak{a}_{B_i}=[N,\frac{-B_i+\sqrt{D}}{2}]=\mathfrak{a}_{B_1}.
$$
Consequently, by Lemma~\ref{lemma_a_b}
\begin{eqnarray*}
J(\mathfrak{a}_1)=J(W_N(\mathfrak{a}_1))
 \Leftrightarrow \mathfrak{a}_{1B} \sim (1) \Leftrightarrow \mathfrak{a}_{iB} \sim (1)  \Leftrightarrow J(\mathfrak{a}_i)=J(W_N(\mathfrak{a}_i)).
\end{eqnarray*}
\end{proof}

\begin{corollary}\label{coro5}
Let $\ell$ be a prime number which splits completely in $K_f$. Let $B$ be an integer such that there exist $u$ and $v$
 satisfied with \eqref{con1}. Further let $\N_B$ be the $N$-system determined by $B$. Then for any $\Fq \in \N_B$,
the polynomial $\Phi_{p_1,p_2}(X,J)$ has a multiple root $\overline{J(\Fa_\Fq)}$ over $\mathbb F_\ell$.
\end{corollary}
\begin{proof}
Our assertion follows from Proposition~\ref{prop2} and Theorem~\ref{theorem4.7}.
\end{proof}

Finally, we give a result for the decomposition of $\Phi_{p_1,p_2}(X,J(\Fa))$ for an ideal $\Fa$ over finite fields.
\begin{proposition} Assume that $p_1$ and $p_2$ satisfy the condition 1 of Corollary~\ref{Cor2.4}. Let $\ell$ be a prime number which splits completely in $K_f$. Let $\Fa$ be an ideal of $\mathcal{O}_f$. Then the polynomial $\Phi_{p_1,p_2}(X,J(\Fa)) \mod \ell$ has at least four linear factors over $\mathbb F_\ell$.
\end{proposition}
\begin{proof}
By Lemmas~\ref{lemma1} and \ref{lemma2}, we know in the ideal class of $\Fa$ there exist four forms $\Fq_B$ of $N$-systems $\N_B$ for four distinct $B \pmod {2N},B^2\equiv D\pmod{4N}~$. By Theorem~\ref{theorem2.3}, we have $\mathfrak{w}_{p_{1},p_{2}}^{s}(\alpha_{\Fq_B})$ are integers of $K_f$. Therefore we have our assertion.
\end{proof}

\section{Example}\label{section5}
We give an example for the result given in Corollary~\ref{coro5}.
Let $D=-56$ and $N=39$. The integer $B=10$ satisfies $B^2\equiv D\pmod{4N}$. Consider the $N$-system $\N_B$. For the integer $B$, there exist $u,v \in \mathbb{Z}$ that satisfy \eqref{con1}. For instance, we can take $u=10,v=1$.

The modular equation $\Phi_{3,13}(X,J)$ and the class polynomial $H_{B,N}(X)$ are given as follows.

{\footnotesize
$
\begin{array}{l}
\Phi_{3,13}(X,J) = 
X^{56}+(704-J)X^{55}+(168568+39J)X^{54}+(14498520-663J)X^{53} \\
\hspace{62pt} +(187807764+6331J)X^{52}+(744637296-35763J)X^{51} \\
\hspace{62pt} +(-6562036+106392J)X^{50}+(-3840625568-18070J)X^{49} \\
\hspace{62pt} +(1058251610-1082016J)X^{48}+(10302034600+3516903J)X^{47} \\
\hspace{62pt} +(4510900472-1278901J)X^{46}+(-34331690432-18277116J)X^{45} \\
\hspace{62pt} +(-7097865034+40532700J)X^{44}+(84188024320+11574823J)X^{43} \\
\hspace{62pt} +(546780176-161476962J)X^{42}+(-154959173464+168751479J)X^{41} \\
\hspace{62pt} +(-12359340101+230086922J)X^{40}+(327081484064-617987682J)X^{39} \\
\hspace{62pt} +(-49301838300+137626281J)X^{38}+(-576339027576+928366231J)X^{37} \\
\hspace{62pt} +(284363953068-959457720J)X^{36}+(735938431592-477589944J)X^{35} \\
\hspace{62pt} +(-558265224452+1429130144J)X^{34}+(-890017323520-466517064J)X^{33} \\
\hspace{62pt} +(977815434427-963208272J)X^{32}+(966995235128+909996295J)X^{31} \\
\hspace{62pt} +(-1755072840368+158515461J)X^{30}+(-345165085024-607329720J)X^{29} \\
\hspace{62pt} +(2218368968890+197238236J)X^{28}+(-911733108784+179445279J)X^{27} \\
\hspace{62pt} +(-1540031876048-140684622J)X^{26}+(1628026178168-6888479J)X^{25} \\
\end{array}
$
}

{\footnotesize
$
\begin{array}{l}
\hspace{62pt} +(261124933147+37909092J)X^{24}+(-1229692547200-8835450J)X^{23} \\
\hspace{62pt} +(462040501468-4070053J)X^{22}+(441029439032+1885689J)X^{21} \\
\hspace{62pt} +(-422841966612+44928J)X^{20}+(-7261052136-111436J)X^{19} \\
\hspace{62pt} +(163453863300+9516J)X^{18}+(-59787354976+740J)X^{17} \\
\hspace{62pt} +(J^2-1486J-26470898021)X^{16}+(24009911816-49J)X^{15} \\
\hspace{62pt} +(-1731574864+29J)X^{14}+(-3926472080+246J)X^{13} \\
\hspace{62pt} +(1333660406-364J)X^{12}+(158103088-221J)X^{11} \\
\hspace{62pt} +(-172600168+650J)X^{10}+(25597000-221J)X^9 \\
\hspace{62pt} +(5195450-364J)X^8+(-2155088+247J)X^7+(177164+26J)X^6 \\
\hspace{62pt} +(39936-52J)X^5+(-9996+13J)X^4+(600-J)X^3+88X^2-16X+1,
\end{array}
$
}

$
H_{B,N}(X)=X^4-2X^3-X^2+2X-1.
$

 It is noted the degree of the $\Phi_{3,13}(X,J)$ in $J$ is 2. We take a prime number $\ell=3593$, which splits completely in $K_1$. Then $H_{B,N}(X)$ decomposes into linear factors over $\mathbb F_\ell$ as follows. 

$H_{B,N}(X)\equiv (X-607)(X-166)(X-3428)(X-2987) \pmod{3593}$.
By substituting the roots of $H_{B,N}(X) \pmod{3593}$ into  $\Phi_{3,13}(X,J)$, we have 
\begin{eqnarray*}
\Phi_{3,13}(607,J) & \equiv & (J-229)^2 \pmod{3593}, \\
\Phi_{3,13}(166,J) & \equiv & (J-2979)^2 \pmod{3593}, \\
\Phi_{3,13}(3428,J) & \equiv & (J-2874)^2 \pmod{3593}, \\
\Phi_{3,13}(2987,J) & \equiv & (J-2696)^2 \pmod{3593}.
\end{eqnarray*}

Let $D(X)$ be the discriminant of $\Phi_{3,13}(X,J)$ as a polynomial in $J$. Then we see $D(X) \equiv 0 \bmod H_{B,N}(X)$. This means that $\Phi_{3,13}(\mathfrak{w}_{p_{1},p_{2}}^{s}(\alpha_{\Fq_i}),J)$ has still a multiple root for every form $\Fq_i\in\N_B$. Thus we also have $J(\mathfrak{a_{\Fq_i}})=J(W_N(\mathfrak{a_{\Fq_i}}))$ for every $i$.


Graduate School of Science \\
Osaka Prefecture University \\
1-1 Gakuen-cho, Naka-ku, Sakai, Osaka 599-8531, Japan \\
dp301004@edu.osakafu-u.ac.jp\vspace{3mm}\\
Graduate School of Science \\
Osaka Prefecture University \\
1-1 Gakuen-cho, Naka-ku, Sakai, Osaka 599-8531, Japan \\
comuta@alg.cias.osakafu-u.ac.jp\vspace{3mm}\\
Faculty of Liberal arts and Sciences \\
Osaka Prefecture University \\
1-1 Gakuen-cho, Naka-ku, Sakai, Osaka 599-8531, Japan \\
ishii@las.osakafu-u.ac.jp\vspace{3mm}\\
\end{document}